\documentclass[twoside,reqno,11pt]{fcaa}

\usepackage{graphicx}
\usepackage{epsfig}
\usepackage{amsthm}
\usepackage{amsmath}
\usepackage{latexsym}
\usepackage{amsfonts}
\usepackage{amssymb}

\newtheoremstyle{theorem}
  {15pt}          
  {15pt}  
  {\sl}  
  {\parindent}
  {\sc}  
  {. }    
  { }    
  {}     
\theoremstyle{theorem}

\newtheoremstyle{defi}
  {15pt}          
  {15pt}  
  {\rm}  
  {\parindent}     
  {\sc}  
  {. }    
  { }    
  {}     
\theoremstyle{defi}

\newtheorem{remark}{Remark}[section]

 \def\proofname{\indent\rm P\ r\ o\ o\ f.}
 \def\proofend{\hfill$\Box$}
 \def\qed{\hfill$\Box$} 

 \usepackage{hyperref} 
 \usepackage{epstopdf}

\usepackage[absolute]{textpos} 
\usepackage{hyperref}

 \def\theequation{\arabic{section}.\arabic{equation}}

 \def\themycopyrightfootnote{\vspace*{3pt}
 \copyright \, Year\,  Diogenes  Co., Sofia
   \par  \noindent pp. xxx--xxx
   \hfill  \vspace*{-36pt}
  }

 \title[Numerical solution of fractional \dots]{Numerical solution of fractional Sturm--Liouville equation in integral form}
 \author[T. Blaszczyk, M. Ciesielski]{Tomasz Blaszczyk $^1$, Mariusz Ciesielski $^2$}

\newtheorem{proposition}{Proposition}[section]

\begin{document}

\begin{textblock}{14}(1.5,1)
\noindent This paper is now published in \\ 
\textbf{Fractional Calculus and Applied Analysis, Vol. 17, No 2 (2014), pp. 307-320}, \\
DOI:10.2478/s13540-014-0170-8, and is available at \url{http://link.springer.com/journal/13540} \\
PDF: \url{http://link.springer.com/article/10.2478/s13540-014-0170-8}
\end{textblock}

 \vbox to 2.5cm { \vfill }


 \bigskip \medskip

 \begin{abstract}

 In this paper a fractional differential equation of the Euler--Lagrange / Sturm--Liouville type is considered. The fractional equation with derivatives of order $\alpha \in \left( 0,1 \right]$ in the finite time interval is transformed to the integral form. Next the numerical scheme is presented. In the final part of this paper examples of numerical solutions of this equation are shown. The convergence of the proposed method on the basis of numerical results is also discussed.

 \medskip


 \smallskip


 \end{abstract}

 \maketitle

 \vspace*{-16pt}



 \section{Introduction}\label{sec:1}

\setcounter{section}{1}
\setcounter{equation}{0}\setcounter{theorem}{0}
The fractional differential equations, both ordinary and partial ones,  are very useful tools for modelling many phenomena in physics, mechanics, control theory, biochemistry, bioengineering and economics~\cite{Blaszczyk05,Hilfer,Leszczynski02,Leszczynski03,Magin,Scalas}. Therefore, the theory of fractional differential equations is an area that has developed extensively over the last decades. In the monographs \cite{Kilbas,Kiryakova,Klimek04,Oldham,Podlubny} one can find a review of methods of solving fractional differential equations. 

In recent years,  subtopic of the theory of fractional differential equations gains importance: it concerns the variational principles for functionals involving fractional derivatives. These principles lead to equations known in the literature as the fractional Euler--Lagrange equations. The equations of this type were derived when fractional integration by parts rule~\cite{Kilbas} has been applied. 

This approach was initiated by Riewe in \cite{Riewe01}, where he used non-integer order derivatives to describe nonconservative systems in mechanics. Next Klimek~\cite{Klimek01} and Agrawal~\cite{Agrawal01} noticed that such equations can be investigated in the sequential approach. A fractional Hamiltonian formalism for the combined fractional calculus of variations was introduced in \cite{ Malinowska03}. In the work \cite{Odzijewicz02} Green theorem for generalized partial fractional derivatives was proved. Other applications of fractional variational principles are presented in~\cite{Agrawal04,Klimek06,Malinowska01,Malinowska02,Odzijewicz01}. 

Recently, the fractional Sturm--Liouville problems were formulated by Klimek and Agrawal in \cite{Klimek07} and Rivero et al. in \cite{Rivero}. Authors in these papers considered several types of the fractional Sturm--Liouville equations and they investigated the eigenvalues and eigenfunctions properties of the fractional Sturm--Liouville operators.

Unfortunately, the fractional Euler--Lagrange / Sturm--Liouville equations contain the composition of the left- and right-sided derivatives. It is an additional drawback for computation of an exact solution (even with simple Lagrangian, see~\cite{Baleanu,Klimek04,Klimek05}). Consequently, numerous studies have been devoted to numerical schemes for the fractional equations (see~\cite{Baleanu1,Baleanu2,Blaszczyk03,Lotfi,Wang}). For numerical methods in the fractional calculus of variations we refer the reader to \cite{Pooseh01,Pooseh02,Pooseh03}.

In our previous works \cite{Blaszczyk01,Blaszczyk03,Blaszczyk04} we proposed numerical scheme on the basis of a finite difference method of solution for a special case of the problem, namely the fractional oscillator equation.  In this paper we propose a numerical solution of the fractional Sturm--Liouville equation.  We investigate a new integral form of this equation and a numerical method of solution of considered equation in conjunction with analysis of a rate of convergence.  Another integral form of the fractional Euler-Lagrange equations (containing the Caputo derivatives) has been recently considered in \cite{Lazo01}.

\section{Statement of the problem and definitions}\label{sec:2}

\setcounter{section}{2}
\setcounter{equation}{0}\setcounter{theorem}{0}
We consider the fractional differential equation with derivatives of order $\alpha \in \left( 0,1 \right]$ in the finite time interval $t \in \left[ 0,b \right]$, for parameter $\lambda \in \mathbb{R}$ and variable potential determined by function $q \left(t\right)$
\begin{equation}
{}^CD_{{b^ - }}^\alpha \,D_{{0^ + }}^\alpha \,f\left( t \right) + \left(\lambda + q \left(t\right) \right) \,f\left( t \right) = 0 ,
\label{eq:fr_osc_diff}
\end{equation}
where $f(t) \in AC[0,b]$ is an unknown function (absolutely continuous on $[0,b]$), satisfying boundary conditions 
\begin{equation}
f\left( 0 \right) = 0,\quad f\left( b \right) = L .
\label{eq:bc}
\end{equation}

According \cite{Kilbas,Oldham,Podlubny} we recall the definitions of the left- and right-sided Riemann-Liouville fractional derivative operators for $\alpha \in \left( 0,1 \right)$  
\begin{eqnarray}
D_{{0^ + }}^\alpha \,f\left( t \right) &:=& D\,I_{{0^ + }}^{1 - \alpha }f\left( t \right)
\label{eq:def_D_RL_left}
\\
D_{{b^ - }}^\alpha \,f\left( t \right) &:=&  - D\,I_{{b^ - }}^{1 - \alpha }f\left( t \right)
\label{eq:def_D_RL_right}
\end{eqnarray}
where $D$ is operator of the first order derivative and operators $I_{{0^ + }}^\alpha$ and $I_{{b^ - }}^\alpha $ are  respectively the left- and right-sided fractional integrals of order $\alpha > 0$ defined by
\begin{eqnarray}
I_{{0^ + }}^\alpha f\left( t \right) &:=& \frac{1}{{\Gamma \left( \alpha  \right)}}\int_0^t {\frac{{f\left( \tau  \right)}}{{{{\left( {t - \tau } \right)}^{1 - \alpha }}}}d\tau } \quad \left( {t > 0} \right)
\label{def_I_left}
\\
I_{{b^ - }}^\alpha f\left( t \right) &:=& \frac{1}{{\Gamma \left( \alpha  \right)}}\int_t^b {\frac{{f\left( \tau  \right)}}{{{{\left( {\tau  - t} \right)}^{1 - \alpha }}}}d\tau } \quad \left( {t < b} \right) ,
\label{def_I_right} 
\end{eqnarray}
whereas operators ${}^CD_{{0^ + }}^\alpha$ and ${}^CD_{{b^ - }}^\alpha$ represent the left- and right-sided Caputo fractional derivatives, respectively. Between both definitions occur the following relationships \cite{Kilbas} (only valid for $\alpha \in (0,1)$) 
\begin{eqnarray}
{}^CD_{{0^ + }}^\alpha \,f\left( t \right) &:=& D_{{0^ + }}^\alpha \,f\left( t \right) - \frac{{{t^{ - \alpha }}}}{{\Gamma \left( {1 - \alpha } \right)}}f\left( 0 \right)
\label{eq:def_D_Caputo_left}
\\
{}^CD_{{b^ - }}^\alpha \,f\left( t \right) &:=& D_{{b^ - }}^\alpha \,f\left( t \right) - \frac{{{{\left( {b - t} \right)}^{ - \alpha }}}}{{\Gamma \left( {1 - \alpha } \right)}}f\left( b \right) .
\label{eq:def_D_Caputo_right}
\end{eqnarray}

In the further part of this paper we will use the following composition rules of  fractional operators (for $\alpha \in \left( 0,1 \right]$) \cite{Kilbas}
\begin{eqnarray}
I_{{0^ + }}^\alpha \,{}^CD_{{0^ + }}^\alpha \,f\left( t \right) &=& f\left( t \right) - f\left( 0 \right) 
\label{eq:composition_rule_left}
\\
I_{{b^ - }}^\alpha \,{}^CD_{{b^ - }}^\alpha \,f\left( t \right) &=& f\left( t \right) - f\left( b \right)
\label{eq:composition_rule_right}
\end{eqnarray}
and the fractional integral of a constant $C$
\begin{equation}
I_{{0^ + }}^\alpha \,C = C\frac{{{t^\alpha }}}{{\Gamma \left( {1 + \alpha } \right)}},\quad I_{{b^ - }}^\alpha C = C\frac{{{{\left( {b - t} \right)}^\alpha }}}{{\Gamma \left( {1 + \alpha } \right)}} .
\label{eq:integration_const}
\end{equation}

In particular, when $\alpha  = 1$, then ${}^CD_{{b^ - }}^1 \,D_{{0^ + }}^1 = -{D^2}$
and for $q \left( t \right) = 0$ Eq. (\ref{eq:fr_osc_diff}) becomes
\begin{equation}
 -{D^2}f\left( t \right) + \lambda \,f\left( t \right) = 0
\end{equation}
and for $\lambda < 0$ (the oscillatory character of solutions of the equation) its analytical solution satisfying  boundary conditions (\ref{eq:bc}) is of the form
\begin{equation}
f\left( t \right) = L\frac{{\sin \left( {\sqrt { - \lambda } \,t} \right)}}{{\sin \left( {\sqrt { - \lambda } \,b} \right)}}, \quad  \lambda \ne  - {\left( {\frac{{k\pi }}{b}} \right)^2}, \quad k \in \mathbb{Z} .
\end{equation}

The formula for the analytical solution of Eq. (\ref{eq:fr_osc_diff}) for $\alpha \in \left( 0,1 \right)$  is rather involved. For example, the analytical solution for $q \left( t \right) = 0$ which contains the series of fractional integrals is presented in \cite{Blaszczyk03,Klimek04,Klimek05}. This makes it difficult to carry out any operations on them. 
There is a problem in calculations of the values of $f$. 
The analytical solution of Eq.~(\ref{eq:fr_osc_diff}) can be expressed by elementary functions only in the special case (see details in further description). 

\section{Transformation of fractional equation to the integral form}\label{sec:3}

\setcounter{section}{3}
\setcounter{equation}{0}\setcounter{theorem}{0}

 \begin{proposition}\label{Pr1}
The equivalent integral form of Eq. (\ref{eq:fr_osc_diff}) with boundary conditions (\ref{eq:bc}) is given as
\begin{equation}
f\left( t \right) + {I_{{0^ + }}^\alpha \,I_{{b^ - }}^\alpha \left(\lambda + q \left(t\right) \right)f\left( t \right) - {{\left( {\frac{t}{b}} \right)}^\alpha }{{\left. {I_{{0^ + }}^\alpha \,I_{{b^ - }}^\alpha \left(\lambda + q \left(t\right) \right)f\left( t \right)} \right|}_{t = b}}}  = {\frac{L}{b^\alpha}}  t^\alpha .
\label{eq:fr_osc_int}
\end{equation}
 \end{proposition}

 \proof 

By using the fractional integral operators  $I_{{0^ + }}^\alpha \,I_{{b^ - }}^\alpha $ acting on Eq. (\ref{eq:fr_osc_diff}), we obtain 
\begin{equation}
I_{{0^ + }}^\alpha \,I_{{b^ - }}^\alpha \,{}^CD_{{b^ - }}^\alpha \,D_{{0^ + }}^\alpha \,f\left( t \right) +  I_{{0^ + }}^\alpha \,I_{{b^ - }}^\alpha \left(\lambda + q \left(t\right) \right)f\left( t \right) = 0 .
\label{eq:intergral_form_1}
\end{equation} 

Next we use the composition rule of operators 
$I_{{b^ - }}^\alpha \,{}^CD_{{b^ - }}^\alpha $ 
(see Eq. (\ref{eq:composition_rule_right})) in Eq. (\ref{eq:intergral_form_1}), thus we get 
\begin{equation}
I_{{0^ + }}^\alpha \,\left( {D_{{0^ + }}^\alpha \,f\left( t \right) - {{\left. {D_{{0^ + }}^\alpha \,f\left( t \right)} \right|}_{t = b}}} \right) + I_{{0^ + }}^\alpha \,I_{{b^ - }}^\alpha \left(\lambda + q \left(t\right) \right)f\left( t \right) = 0
\label{eq:intergral_form_2}
\end{equation}
The above equation contains an unknown value ${\left. {D_{{0^ + }}^\alpha \,f\left( t \right)} \right|}_{t = b}$ and for the unknown function $f(t)$ this value is treated here as a constant. 

Using again the composition rule of operators 
$I_{{0^ + }}^\alpha \, D_{{0^ + }}^\alpha $ 
(see Eq. (\ref{eq:composition_rule_left}) and the fact that if $f(0) = 0$, then 
$ D_{{0^ + }}^\alpha \,f\left( t \right) = {}^CD_{{0^ + }}^\alpha \,f\left( t \right) $ 
in Eq. (\ref{eq:def_D_Caputo_left}), hence $I_{{0^ + }}^\alpha \,D_{{0^ + }}^\alpha \,f\left( t \right) = f\left( t \right)$) 
and the fractional integral of a constant (\ref{eq:integration_const}), we obtain the following form of equation
\begin{equation}
f\left( t \right) - {\left. {D_{{0^ + }}^\alpha \,f\left( t \right)} \right|_{t = b}}\frac{{{t^\alpha }}}{{\Gamma \left( {\alpha  + 1} \right)}} + I_{{0^ + }}^\alpha \,I_{{b^ - }}^\alpha \left(\lambda + q \left(t\right) \right)f\left( t \right) = 0.
\label{eq:intergral_form_3}
\end{equation}

In order to determine the value ${\left. {D_{{0^ + }}^\alpha \,f\left( t \right)} \right|}_{t = b}$ we substitute the value $t=b$ into Eq. (\ref{eq:intergral_form_3})
\begin{equation}
f\left( b \right) - {\left. {D_{{0^ + }}^\alpha \,f\left( t \right)} \right|_{t = b}}\frac{{{b^\alpha }}}{{\Gamma \left( {\alpha  + 1} \right)}} + {\left. {I_{{0^ + }}^\alpha \,I_{{b^ - }}^\alpha \left(\lambda + q \left(t\right) \right)f\left( t \right)} \right|_{t = b}} = 0
\label{eq:intergral_form_4}
\end{equation}
and obtain 
\begin{equation}
{\left. {D_{{0^ + }}^\alpha \,f\left( t \right)} \right|_{t = b}} = \frac{{\Gamma \left( {\alpha  + 1} \right)}}{{{b^\alpha }}}\left( {f\left( b \right) + {{\left. {I_{{0^ + }}^\alpha \,I_{{b^ - }}^\alpha \left(\lambda + q \left(t\right) \right)f\left( t \right)} \right|}_{t = b}}} \right).
\label{eq:intergral_form_5}
\end{equation}
Next we substitute the right-hand side of (\ref{eq:intergral_form_5}) into Eq. (\ref{eq:intergral_form_3}) and get the integral form of Eq. (\ref{eq:fr_osc_diff})
\begin{align}
f\left( t \right) &+  I_{{0^ + }}^\alpha \,I_{{b^ - }}^\alpha \left(\lambda + q \left(t\right) \right)f\left( t \right) - {{\left( {\frac{t}{b}} \right)}^\alpha }{{\left. {I_{{0^ + }}^\alpha \,I_{{b^ - }}^\alpha \left(\lambda + q \left(t\right) \right)f\left( t \right)} \right|}_{t = b}} \nonumber \\
&  =  {\left( {\frac{t}{b}} \right)^\alpha } {f\left( b \right) }  .
\label{eq:fr_osc_int_general}
\end{align}
Taking into account the boundary conditions (\ref{eq:bc}) we obtain Eq. (\ref{eq:fr_osc_int}).
 \proofend 

 \begin{remark}\label{Re1}
One can easily check that ${\left. {I_{{0^ + }}^\alpha \,I_{{b^ - }}^\alpha \left(\lambda + q \left(t\right) \right)f\left( t \right)} \right|}_{t = 0} = 0$ (on the basis of definitions (\ref{def_I_left}) and (\ref{def_I_right})). If we put $t=0$ and $t=b$ into Eq. (\ref{eq:fr_osc_int}) we can confirm that this equation fulfils the boundary conditions (\ref{eq:bc}). In particular, for $\lambda + q \left( t \right) = 0$ Eq. (\ref{eq:fr_osc_int}) simplifies to the form
\begin{equation}
f\left( t \right) =  {\frac{L}{b^\alpha}}  t^\alpha .
\label{eq:fr_osc_int_lam0}
\end{equation}
\end{remark}

\section{Numerical algorithm}\label{sec:4}

\setcounter{section}{4}
\setcounter{equation}{0}\setcounter{theorem}{0}
In this section we present a numerical scheme  for Eq. (\ref{eq:fr_osc_int}). We introduce the homogeneous grid of nodes (with the constant time step $\Delta t = b/n$, where $n+1$ is a number of nodes): $0 = t_0 < t_1 < ... < t_i < t_{i+1} < ... < t_n = b$, and $t_i = i \, \Delta t$. In order to simplify notation we will introduce a new function $\phi \left( t \right) \equiv \left(\lambda + q \left(t\right) \right)f\left( t \right)$ and we denote the values of functions $f(t)$, $q(t)$ and $\phi (t)$ at the node $t_i$ by $f_i = f(t_i)$, $q_i = q(t_i)$ and $ {\phi }_i = \phi (t_i) = (\lambda + q_i)f_i$.

Now we determine the numerical schemes  of integration \cite{Diethelm01,Oldham, Press} for both fractional integral operators occurring in Eq. (\ref{eq:fr_osc_int}).

At node $t_0$  we have
${\left. {I_{{0^ + }}^\alpha \phi\left( t \right)} \right|_{t = {t_0}}} = 0$.
Discrete form of the integral operator (\ref{def_I_left}) at nodes $t_i$ for $i = 1, 2, ..., n$ is approximated by the formula
\begin{eqnarray}
{\left. {I_{{0^ + }}^\alpha \phi\left( t \right)} \right|_{t = {t_i}}} 
& = & \frac{1}{{\Gamma \left( \alpha  \right)}}\int_0^{{t_i}} {\frac{{\phi\left( \tau  \right)}}{{{{\left( {{t_i} - \tau } \right)}^{1 - \alpha }}}}d\tau }  = \frac{1}{{\Gamma \left( \alpha  \right)}}\sum\limits_{j = 0}^{i - 1} {\int_{{t_j}}^{{t_{j + 1}}} {\frac{{\phi\left( \tau  \right)}}{{{{\left( {{t_i} - \tau } \right)}^{1 - \alpha }}}}d\tau } } \nonumber \\
& \approx & \frac{1}{{\Gamma \left( \alpha  \right)}}\sum\limits_{j = 0}^{i - 1} {\frac{{{\phi_j} + {\phi_{j + 1}}}}{2}\int_{j\,\Delta t}^{\left( {j + 1} \right)\Delta t} {\frac{1}{{{{\left( {i\,\Delta t - \tau } \right)}^{1 - \alpha }}}}d\tau } } \nonumber \\
& = & \frac{{{{\left( {\Delta t} \right)}^\alpha }}}{{2\Gamma \left( {\alpha  + 1} \right)}}\sum\limits_{j = 0}^{i - 1} {\left( {{\phi_j} + {\phi_{j + 1}}} \right)\left( {{{\left( {i - j} \right)}^\alpha } - {{\left( {i - j - 1} \right)}^\alpha }} \right)} \nonumber \\
& = & \sum\limits_{j = 0}^i {{\phi_j} \,{w_{i,j}}} ,
\end{eqnarray}
where the coefficients $w_{i,j}$ (also including the case for $i=0$) are as follows
\begin{equation}
{w_{i,j}} = \frac{{{{\left( {\Delta t} \right)}^\alpha }}}{{2\Gamma \left( {\alpha  + 1} \right)}}\left\{ {\begin{array}{*{20}{l}}
0&{{\mbox{for }}i = 0{\mbox{ and }}j = 0}\\
{{i^\alpha } - {{\left( {i - 1} \right)}^\alpha }}&{{\mbox{for }}i > 0{\mbox{ and }}j = 0}\\
{{{\left( {i - j + 1} \right)}^\alpha } - {{\left( {i - j - 1} \right)}^\alpha }}&{{\mbox{for }}i > 0{\mbox{ and }} 0 < j < i }\\
1&{{\mbox{for }}i > 0{\mbox{ and }}j = i} 
\end{array}} \right. .
\end{equation}

We determine discrete form of the fractional integral operator (\ref{def_I_right}) in a similar way. This operator at node $t_n$ is equal to ${\left. {I_{{b^ - }}^\alpha \phi\left( t \right)} \right|_{t = {t_n}}} = 0$, whereas at nodes $t_i$, $i = 0,1,...,n-1$, the discrete values are determined by the formula
\begin{eqnarray}
{\left. {I_{{b^ - }}^\alpha \phi\left( t \right)} \right|_{t = {t_i}}} 
&=& \frac{1}{{\Gamma \left( \alpha  \right)}}\int_{{t_i}}^{{t_n}} {\frac{{\phi\left( \tau  \right)}}{{{{\left( {\tau  - {t_i}} \right)}^{1 - \alpha }}}}d\tau }  = \frac{1}{{\Gamma \left( \alpha  \right)}}\sum\limits_{j = i}^{n - 1} {\int_{{t_j}}^{{t_{j + 1}}} {\frac{{\phi\left( \tau  \right)}}{{{{\left( {\tau  - {t_i}} \right)}^{1 - \alpha }}}}d\tau } } \nonumber \\
&\approx & \frac{1}{{\Gamma \left( \alpha  \right)}}\sum\limits_{j = i}^{n - 1} {\frac{{{\phi_j} + {\phi_{j + 1}}}}{2}\int_{j\,\Delta t}^{\left( {j + 1} \right)\Delta t} {\frac{1}{{{{\left( {\tau  - i\,\Delta t} \right)}^{1 - \alpha }}}}d\tau } } \nonumber \\
&=& \frac{{{{\left( {\Delta t} \right)}^\alpha }}}{{2\Gamma \left( {\alpha  + 1} \right)}}\sum\limits_{j = i}^{n - 1} {\left( {{\phi_j} + {\phi_{j + 1}}} \right)\left( {{{\left( {j - i + 1} \right)}^\alpha } - {{\left( {j - i} \right)}^\alpha }} \right)} \nonumber \\
&=& \sum\limits_{j = i}^n {{\phi_j} \,{v_{i,j}}} ,
\end{eqnarray}
where the coefficients $v_{i,j}$ (together with the case for $i=n$) have the form
\begin{equation}
{v_{i,j}} = \frac{{{{\left( {\Delta t} \right)}^\alpha }}}{{2\Gamma \left( {\alpha  + 1} \right)}}\left\{ {\begin{array}{*{20}{l}}
0&{{\mbox{for }}i = n{\mbox{ and }}j = n}\\
{{{\left( {n - i} \right)}^\alpha } - {{\left( {n - i - 1} \right)}^\alpha }}&{{\mbox{for }}i < n{\mbox{ and }}j = n}\\
{{{\left( {j - i + 1} \right)}^\alpha } - {{\left( {j - i - 1} \right)}^\alpha }}&{{\mbox{for }}i < n{\mbox{ and }} i < j < n}\\
1&{{\mbox{for }}i < n{\mbox{ and }}j = i}
\end{array}} \right. .
\end{equation}
The discrete form of the composition of both operators $I_{{0^ + }}^\alpha \,I_{{b^ - }}^\alpha \phi\left( t \right)$ at nodes $t=t_i$ for $i = 0,1,...,n$, has the following form
\begin{equation}
{\left. {I_{{0^ + }}^\alpha \,I_{{b^ - }}^\alpha \phi\left( t \right)} \right|_{t = {t_i}}} \approx \sum\limits_{j = 0}^i {{w_{i,j}}\sum\limits_{k = j}^n {{\phi_k} \,{v_{j,k}}} } ,
\end{equation}
or
\begin{equation}
{\left. {I_{{0^ + }}^\alpha \,I_{{b^ - }}^\alpha \left(\lambda + q \left(t\right) \right)f\left( t \right)} \right|_{t = {t_i}}} \approx \sum\limits_{j = 0}^i {{w_{i,j}}\sum\limits_{k = j}^n {{\left(\lambda + q_k\right) f_k} \,{v_{j,k}}} } .
\end{equation}
One can note that at the node $t_0$ we have ${\left. {I_{{0^ + }}^\alpha \,I_{{b^ - }}^\alpha \phi\left( t \right)} \right|_{t = {t_0}}} = 0$.

Now we present the discrete form of the integral equation (\ref{eq:fr_osc_int}).
The solution one can write in the form of the system of $n+1$ linear equations.
For every grid node $t_i$, $i = 0,1,...,n$, we write the following equation
\begin{align}
{f_i} &+ {\sum\limits_{j = 0}^i {{w_{i,j}}\sum\limits_{k = j}^n {{\left(\lambda + q_k\right) f_k}{v_{j,k}}} }  - {{\left( {\frac{i}{n}} \right)}^\alpha }\sum\limits_{j = 0}^n {{w_{n,j}}\sum\limits_{k = j}^n {{\left(\lambda + q_k\right) f_k}{v_{j,k}}} } } \nonumber \\
& = {\left( {\frac{i}{n}} \right)^\alpha }L .
\label{eq:disc_form}
\end{align}
Analysing the above system of equations, created equations for node indexes $i=0$ and $i=n$ one can reduce to the forms $f_0=0$ and $f_n=L$, respectively. In this way the obtained system of linear equations   can be solved numerically.

\section{Simulation results and numerical error analysis}\label{sec:5}

\setcounter{section}{5}
\setcounter{equation}{0}\setcounter{theorem}{0}
In this section we present the results of calculation obtained by our numerical approach to the fractional Sturm--Liouville equation (\ref{eq:fr_osc_diff}).
In order to numerically solve the system of equations (\ref{eq:disc_form}) we used the LUP decomposition method \cite{Press}.
We present several examples of calculations for different values of parameters $\alpha$, $\lambda$ and forms of function $q(t)$. In all these examples we assumed: $b=1$ and $L=1$ in the right boundary condition. For all presented graphs of functions (Figures 1-3), in the calculations we assume that the time domain $t \in [0,1]$ has been divided into $n=2048$ subintervals.

\subsection{Example of results}
Figures 1 and 2 show the example graphs of solutions  of Eq. (\ref{eq:fr_osc_diff}) for $q(t) = 0$ (the case of the fractional oscillator equation). 
Figure 1 presents solutions for $\lambda \in \{-3, -10, -20, -25\}$  and variable values of parameter $\alpha$.  
We can see  how changes of parameter $\alpha$ influence on the frequency of oscillations in comparison to the classical oscillator equation ($\alpha =1$). 
Whereas in the Figure 3 we show the influence of parameter $\lambda \in \{-5, -7.5, -10, -12.5\}$ at the constant values of $\alpha = 0.6$ and $\alpha = 0.8$ on the solution.

\begin{figure}[!h]
\centering
\resizebox*{0.495\textwidth}{!}{\includegraphics{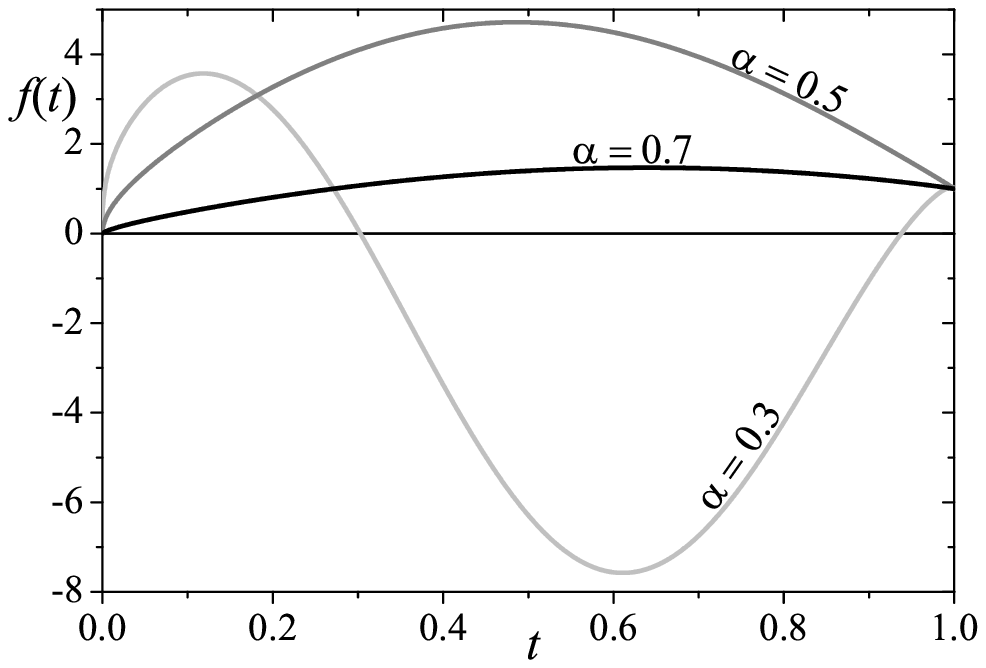}}
\resizebox*{0.495\textwidth}{!}{\includegraphics{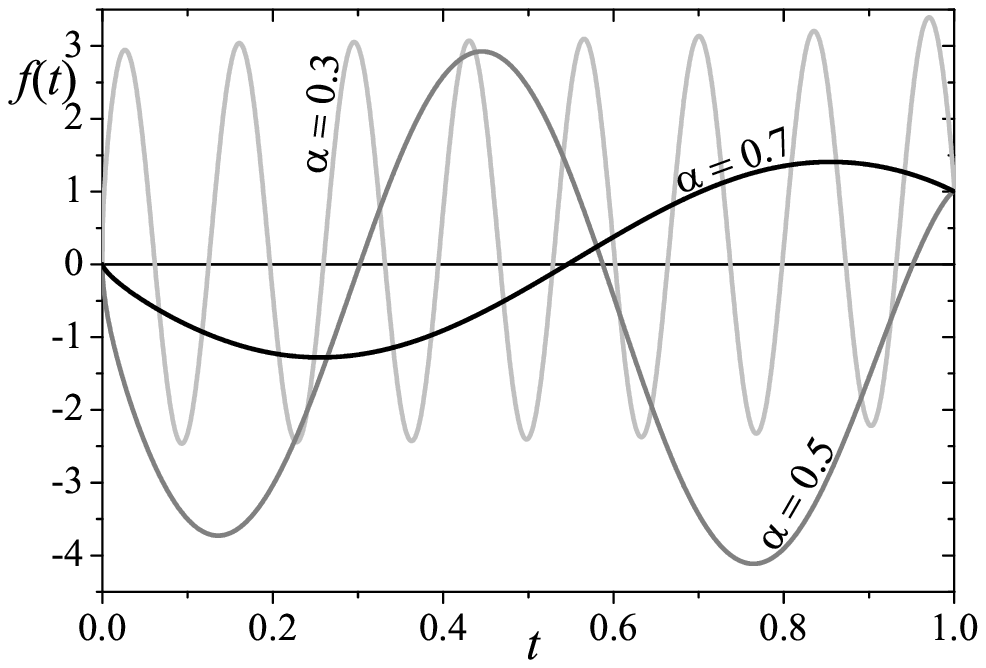}}

\resizebox*{0.495\textwidth}{!}{\includegraphics{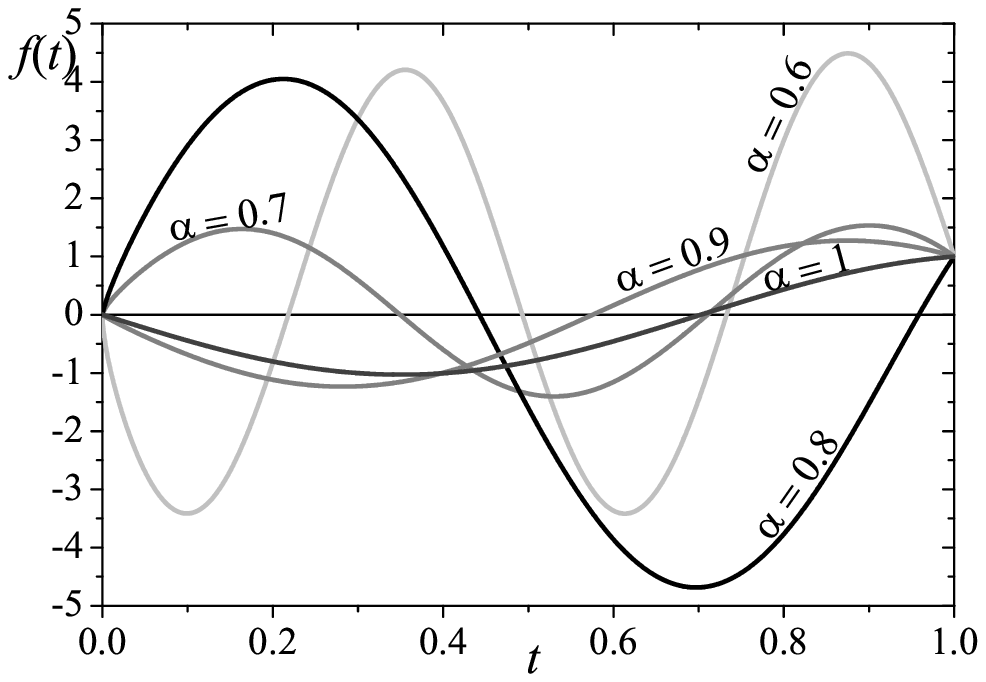}}
\resizebox*{0.495\textwidth}{!}{\includegraphics{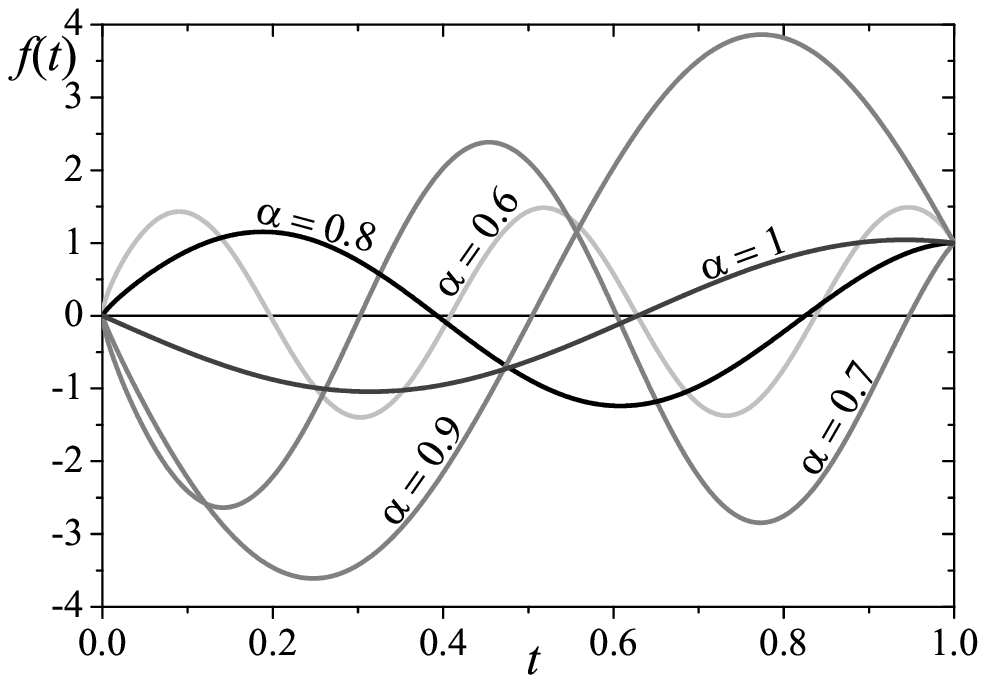}}
\label{fig1}
\vspace{-0.3cm}
\caption{Numerical solution of Eq. (\ref{eq:fr_osc_diff}) for different parameters $\alpha $, $q(t) = 0$,  $b = 1$, $L = 1$ and 
$\lambda = -3$ (left/top), $\lambda = -10$ (right/top), $\lambda = -20$ (left/bottom) and $\lambda = -25$ (right/bottom)}
\end{figure}

\begin{figure}[!h]
\centering
\resizebox*{0.495\textwidth}{!}{\includegraphics{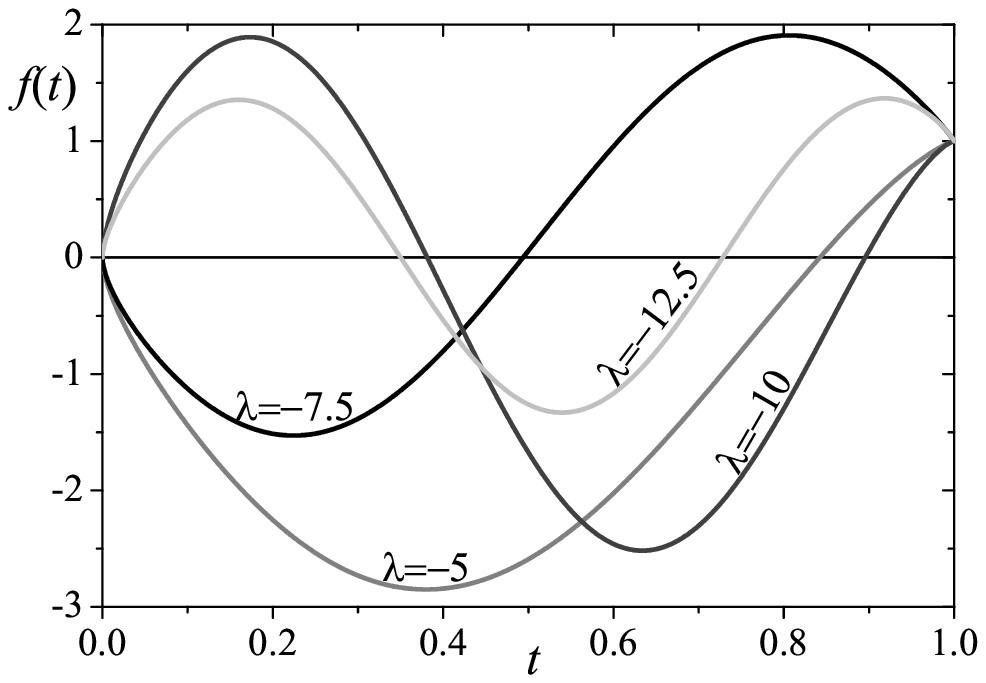}}
\resizebox*{0.495\textwidth}{!}{\includegraphics{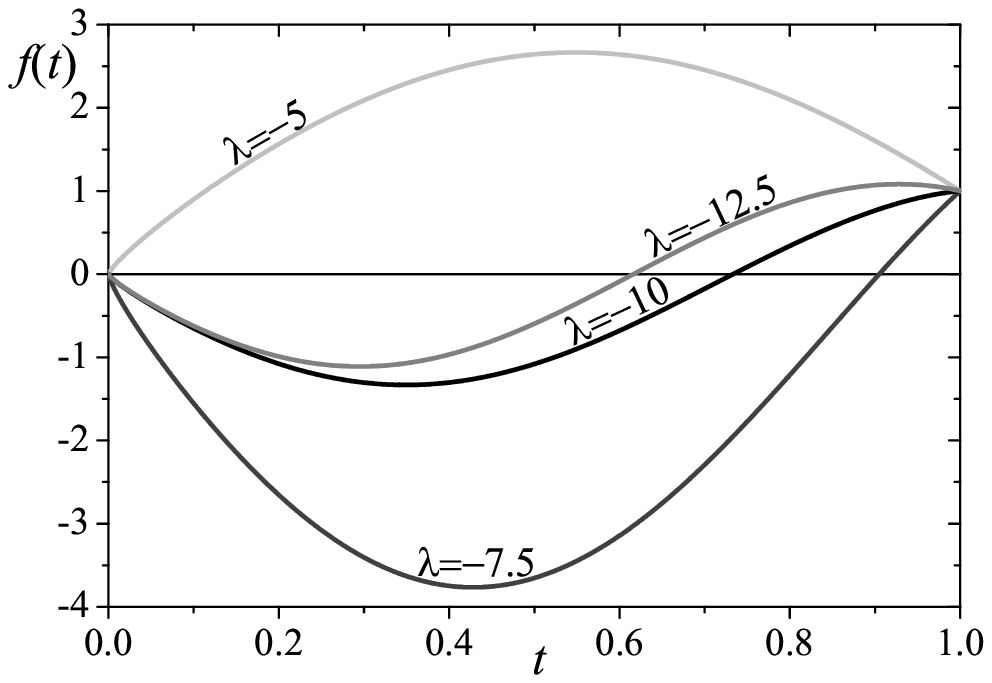}}
\label{fig2}
\vspace{-0.3cm}
\caption{Numerical solution of Eq. (\ref{eq:fr_osc_diff}) for $\lambda \in \{-5, \allowbreak -7.5,  -10, -12.5\}$,  $b = 1$, $L = 1$, $q(t)=0$ and $\alpha = 0.6$  (left-side),  $\alpha = 0.8$ (right-side)}
\end{figure}

In the last Figure 3 we present solutions of Eq. (\ref{eq:fr_osc_diff}) for different form of function $q(t) \neq 0$ (the case of  the fractional Sturm--Liouville equation).
We show the influence of parameters $\alpha$, $\lambda$ and forms of function $q(t)$ on the solution. 
\begin{figure}[!h]
\centering
\resizebox*{0.495\textwidth}{!}{\includegraphics{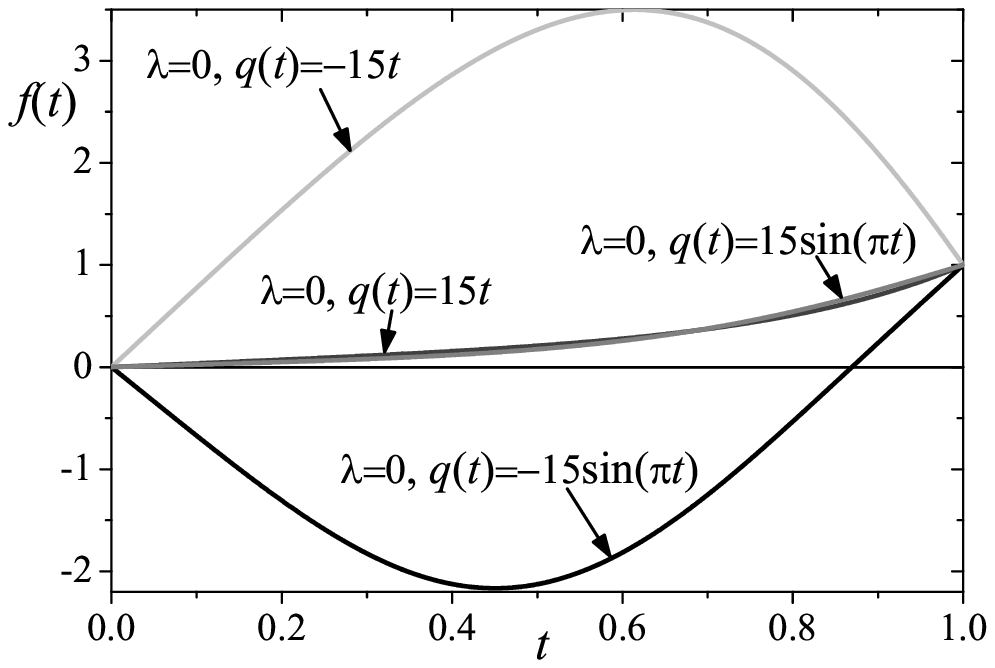}}
\resizebox*{0.495\textwidth}{!}{\includegraphics{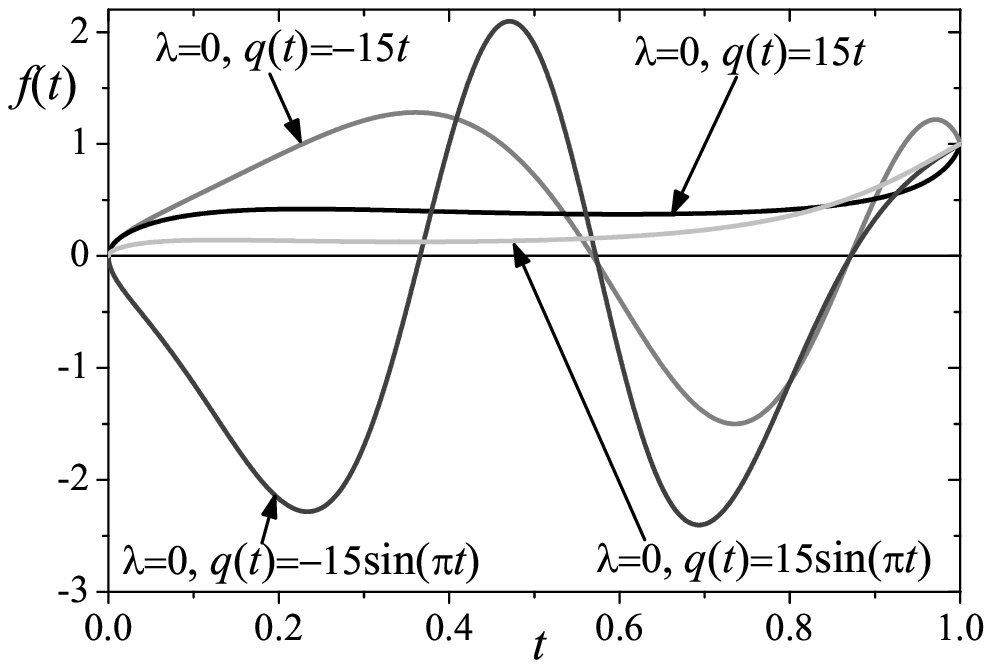}}
\resizebox*{0.495\textwidth}{!}{\includegraphics{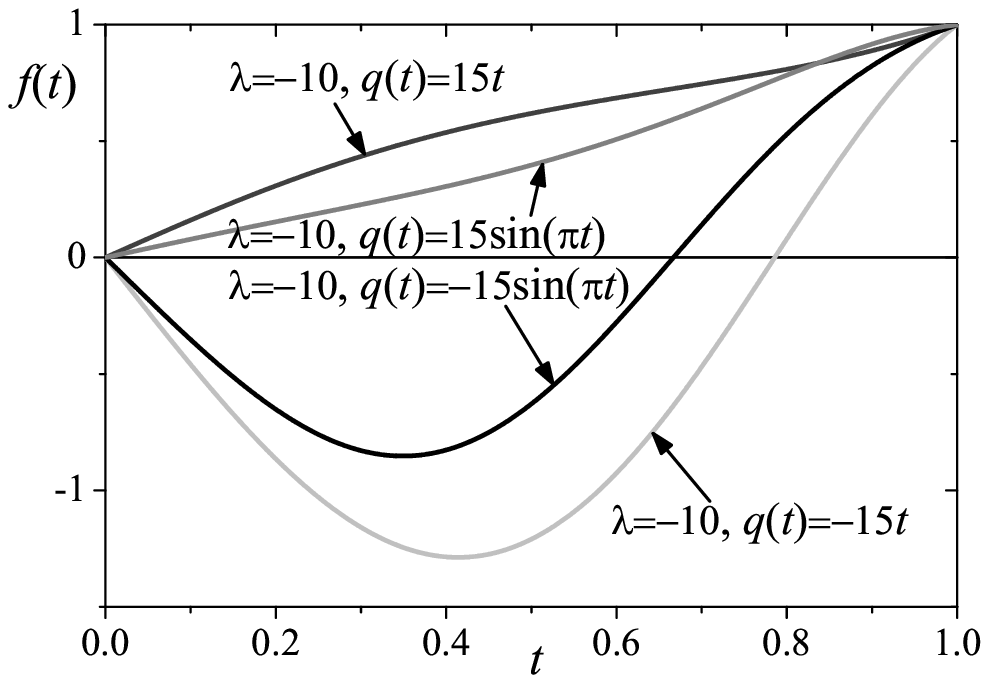}}
\resizebox*{0.495\textwidth}{!}{\includegraphics{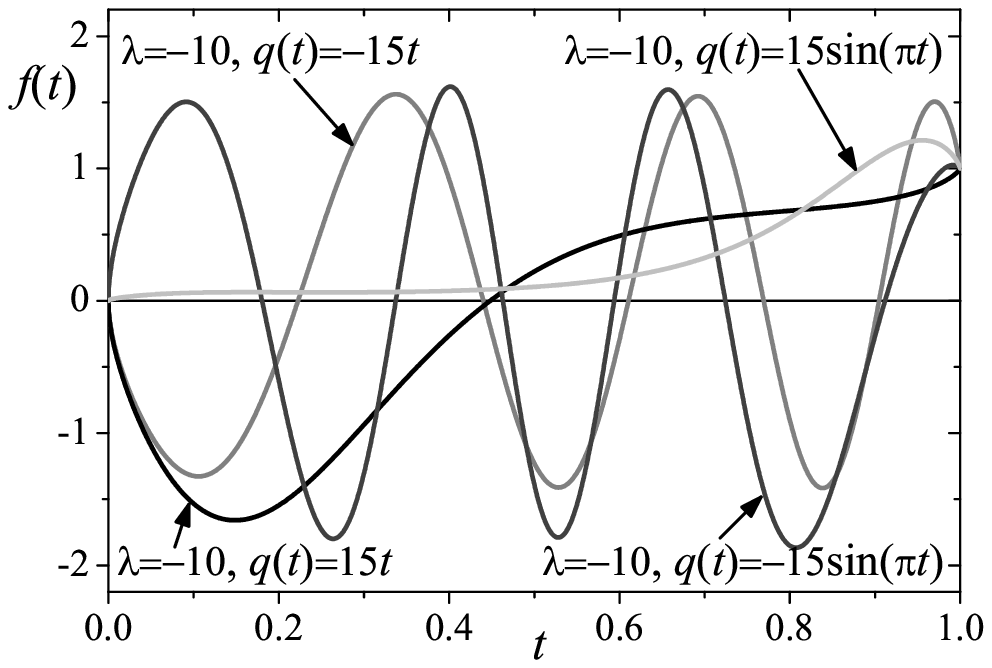}}
\label{fig3}
\vspace{-0.3cm}
\caption{Numerical solution of Eq. (\ref{eq:fr_osc_diff}) for different parameters $\lambda $, forms of function $q(t)$ and $\alpha = 1$ (left-side),  $\alpha = 0.5$ (right-side)}
\end{figure}

\subsection{Error analysis}
Next we analyse errors and convergence of the numerical scheme (\ref{eq:disc_form}) for $q(t)=0$, any $\lambda$ and $ \alpha \in \left(0, 1 \right]$.
When analytical solution is not available, the rate of convergence $p = p_i(\Delta t,\alpha ,\lambda)$ at nodes $t_i$, for fixed parameters $\alpha$, $\lambda$ and variable values of $\Delta t$, can be determined from the following formula (see proposition~\cite{Agrawal05})
\begin{equation}
R_i^{\left( {\Delta t,\alpha ,\lambda } \right)} = \frac{{f_i^{\left( {\Delta t,\alpha ,\lambda } \right)} - f_i^{\left( {2\Delta t,\alpha ,\lambda } \right)}}}{{f_i^{\left( {\Delta t/2,\alpha ,\lambda } \right)} - f_i^{\left( {\Delta t,\alpha ,\lambda } \right)}}} = {2^{p_i\left( {\Delta t,\alpha ,\lambda } \right)}} .
\end{equation}
We thus have
\begin{equation}
p_i\left( {\Delta t,\alpha ,\lambda } \right) = {\log _2}\frac{{f_i^{\left( {\Delta t,\alpha ,\lambda } \right)} - f_i^{\left( {2\Delta t,\alpha ,\lambda } \right)}}}{{f_i^{\left( {\Delta t/2,\alpha ,\lambda } \right)} - f_i^{\left( {\Delta t,\alpha ,\lambda } \right)}}} .
\end{equation}

We present numerical values at three selected nodes and rate of convergence for $ \alpha \in \{ 0.3, 0.5, 0.7 \}$ and $\lambda = -3$, $q(t)=0$ in Table 1. In Table 2 the numerical values at three selected nodes and rates of convergence for $ \alpha = 0.6 $, $\lambda \in \{-5, -7.5, -10 \}$ and $q(t)=0$ are shown. The values from Tables 1 and 2  are also shown in plots - Figures 1 and 2, respectively.

\begin{table}[h]
\centering
\caption{Numerical values of $f$ at nodes $t_i$, $i \in \{n/4, n/2, \allowbreak 3n/4\}$ 
and rates of convergence $p$ for parameters $\lambda = -3$, $q(t)=0$, $b=1$, $L=1$}
\label{tab01}
\small
\begin{tabular}{|l|l|l|l|l|l|l|l|}
\hline
~& ~& \multicolumn{2}{|c|}{$t=0.25$}&\multicolumn{2}{|c|}{$t=0.5$}& \multicolumn{2}{|c|}{$t=0.75$}\\ 
\hline
$\alpha$& $\Delta t = 1/n$& $f_{n/4}$ & $p$ & $f_{n/2}$ & $p$  & $f_{3n/4}$ & $p$ \\
\hline
0.3             & 1/256 & 1.53966755 & - & -6.00003222 & - & -5.30712828 & - \\ 
~              & 1/512 & 1.58297521 & 1.08 & -6.17593877 & 1.12 & -5.50843511 & 1.12 \\ 
~   & 1/1024 & 1.60340765 & 1.18 & -6.25702726 & 1.19 & -5.60111461 & 1.20 \\ 
~   & 1/2048 & 1.61245514 & 1.22 & -6.29245343 & 1.24 & -5.64157016 & 1.24 \\ 
~   & 1/4096 & 1.61632922 & 1.25 & -6.30749294 & 1.26 & -5.65873385 & 1.26 \\ 
~   & 1/8192 & 1.61795763 & - & -6.31377732 & - & -5.66590217 & - \\ 
\hline
0.5            & 1/256 & 3.73516937 & - & 4.73052344 & - & 3.57465003 & - \\ 
~              & 1/512 & 3.72842313 & 1.38& 4.72211467 & 1.39& 3.56919099 & 1.39\\ 
~   & 1/1024 & 3.72583618 & 1.42& 4.71890721 & 1.42& 3.56710826 & 1.42\\ 
~   & 1/2048 & 3.72486755 & 1.44& 4.71771054 & 1.44& 3.56633103 & 1.44\\ 
~   & 1/4096 & 3.72451075 & 1.46& 4.71727084 & 1.46& 3.56604538 & 1.46\\ 
~   & 1/8192 & 3.72438081 & - & 4.71711100 & - & 3.56594153 & - \\ 
\hline
0.7            & 1/256 & 0.94678077 & - & 1.40241981 & - & 1.42764144 & - \\ 
~              & 1/512 & 0.94671534 & 1.52& 1.40231616 & 1.55& 1.42754955 & 1.57\\ 
~   & 1/1024 & 0.94669252 & 1.57& 1.40228081 & 1.59& 1.42751850 & 1.60 \\ 
~   & 1/2048 & 0.94668482 & 1.60& 1.40226907 & 1.61& 1.42750826 & 1.62\\ 
~   & 1/4096 & 0.94668228 & 1.63& 1.40226523 & 1.64& 1.42750493 & 1.64\\ 
~   & 1/8192 & 0.94668146 & - & 1.40226400 & - & 1.42750386 & - \\ 
\hline
\end{tabular}
\end{table}

\begin{table}[h]
\centering
\caption{Numerical values of $f$ at nodes $t_i$, $i \in \{n/4, n/2, \allowbreak 3n/4\}$
and rates of convergence $p$ for parameters $\alpha = 0.6$, $q(t)=0$, $b=1$, $L=1$}\label{tab02}
\small
\begin{tabular}{|l|l|l|l|l|l|l|l|}
\hline
~& ~& \multicolumn{2}{|c|}{$t=0.25$}&\multicolumn{2}{|c|}{$t=0.5$}& \multicolumn{2}{|c|}{$t=0.75$}\\ 
\hline
$\lambda$ & $\Delta t = 1/n$& $f_{n/4}$ & $p$ & $f_{n/2}$ & $p$  & $f_{3n/4}$ & $p$ \\ 
\hline
-5              & 1/256  & -2.16736188& - & -2.58746851 & - & -0.79882549 & - \\ 
~              & 1/512  & -2.16934247& 1.46 & -2.59034981 & 1.44 & -0.80084057 & 1.45 \\ 
~           & 1/1024  & -2.17006290& 1.50 & -2.59140830 & 1.49 & -0.80157662 & 1.49 \\ 
~          & 1/2048  & -2.17031750& 1.53 & -2.59178487 & 1.52 & -0.80183747 & 1.52 \\ 
~          & 1/4096  & -2.17040578& 1.55 & -2.59191604 & 1.54 & -0.80192808 & 1.54 \\ 
~          & 1/8192  & -2.17043598& - & -2.59196106 & - & -0.80195912 & - \\ 
\hline
-7.5           & 1/256 & -1.51542247& - & 0.06057150 & - & 1.83620282 & - \\ 
~               & 1/512 & -1.51337627& 1.46& 0.05880539 & 1.40& 1.83222891 & 1.45\\ 
~   & 1/1024 & -1.51263480& 1.50& 0.05815391 & 1.46& 1.83078071 & 1.50\\ 
~   & 1/2048 & -1.51237316& 1.53& 0.05792141 & 1.50& 1.83026791 & 1.53\\ 
~   & 1/4096 & -1.51228251& 1.54& 0.05784024 & 1.53& 1.83008983 & 1.55\\ 
~   & 1/8192 & -1.51225151& - & 0.05781234 & - & 1.83002882 & - \\ 
\hline
-10            & 1/256 & 1.58290914 & - & -1.66780189 & - & -1.87108209 & - \\ 
~               & 1/512 & 1.58736775 & 1.42& -1.67135428 & 1.46& -1.87951726 & 1.42\\ 
~   & 1/1024 & 1.58904541 & 1.48& -1.67267356 & 1.51& -1.88265567& 1.48\\ 
~   & 1/2048 & 1.58965039 & 1.51& -1.67314544 & 1.53& -1.88377958 & 1.52\\ 
~   & 1/4096 & 1.58986289 & 1.54& -1.67331027 & 1.56& -1.88417251 & 1.54\\ 
~   & 1/8192 & 1.58993623 & - & -1.67336695 & - & -1.88430769 & - \\ 
\hline
\end{tabular}
\end{table}

Analysing the values in the above tables, one can observe that the rate of convergence $p$ 
is dependent on the fractional order $\alpha$ and does not depend on parameter 
$\lambda$. The rate of convergence $p$ is close to $1 + \alpha$.

\section{Conclusions}\label{sec:6}

\setcounter{section}{6}
\setcounter{equation}{0}\setcounter{theorem}{0}
In this paper the fractional Sturm--Liouville equation with derivatives of order 
$\alpha \in \left( 0,1 \right]$ in the finite time interval 
$t \in \left[ 0,b \right]$ is considered. This equation was
transformed using the composition rules for 
fractional integrals and derivatives to the integral form. 
Next the discrete form of the 
integral equation was presented as the system of linear algebraic equations.
The obtained system of equations was solved numerically. 
The equation was solved for derivatives of different orders $\alpha$,  
different values of parameter $\lambda$ and different forms of function $q(t)$. 
The presented results showed the influence these values on the character (i.e. the occurrence of oscillations) of the solution. 
One can note that for $q(t)=0$ the 
oscillations occur only for $\lambda < 0$.
The number of oscillations increases when the value of order $\alpha$ 
decreases for fixed value of parameter $\lambda$. Similarly, for 
fixed order $\alpha$, the number of oscillations increases when 
the value of parameter $\lambda$ decreases. In order to ensure 
stability of the computation, the convergence study of the numerical
scheme was conducted. The rate of convergence $p$ was estimated to be 
close to $1 + \alpha$. 

\smallskip
\section*{Acknowledgements}

This work was supported by the Czestochowa University of Technology Grant 
Number BS/MN 1-105-302/13/P.




 \bigskip \smallskip

 \it

 \noindent
$^1$ Institute of Mathematics \\
Czestochowa University of Technology \\
al. Armii Krajowej 21, 42-201 Czestochowa, POLAND  \\[4pt]
e-mail: tomasz.blaszczyk@im.pcz.pl
\hfill Received: June 18, 2013 \\[12pt]

\noindent
$^2$ Institute of Computer and Information Sciences\\
Czestochowa University of Technology\\
ul. Dabrowskiego 73,  42-201 Czestochowa, POLAND \\[4pt]
e-mail: mariusz.ciesielski@icis.pcz.pl


\begin{thebibliography}{99}
 \normalsize


\bibitem{Agrawal01} O. P. Agrawal, Formulation of Euler-Lagrange equations for
                    fractional variational problems,
                    \emph{J. Math. Anal. Appl.} \textbf{272} (2002), 368-379.
\bibitem{Agrawal04} O. P. Agrawal, S. I. Muslih and D. Baleanu, 
                    Generalized variational calculus in terms of 
                    multi-parameters fractional derivatives,
                    \emph{Commun. Nonlinear Sci. Numer. Simulat.} \textbf{16} (2011), 4756-4767.
\bibitem{Agrawal05} O. P. Agrawal, M. M. Hasan and X. W. Tangpong, 
                    A numerical scheme for a class of parametric problem of fractional 
                    variational calculus,
                    \emph{J. Comput. Nonlinear Dyn.} \textbf{7} (2012), 021005-1--021005-6.
\bibitem{Baleanu1} D. Baleanu, K. Diethelm, E. Scalas, J.J. Trujillo, 
                   {\it  Fractional Calculus Models and Numerical Methods}, World Scientific, Singapore (2012).
\bibitem{Baleanu} D. Baleanu and J. J. Trujillo, On exact solutions of a class
                    of fractional Euler-Lagrange equations,
                    \emph{Nonlinear Dyn.} \textbf{52} (2008), 331-335.
\bibitem{Baleanu2}  D. Baleanu, I. Petras, J. H. Asad and M. P. Velasco, 
                    Fractional Pais-Uhlenbeck Oscillator,
                    \emph{International Journal of Theoretical Physics} 
                    \textbf{51(4)} (2012), 1253-1258.
\bibitem{Blaszczyk01} T. Blaszczyk and M. Ciesielski, Fractional Euler-Lagrange 
                    equations -- numerical solutions and applications of 
                    reflection operator, \emph{Scientific Research of the Institute 
                    of Mathematics and Computer Science} \textbf{2(9)} (2010), 17-24.
\bibitem{Blaszczyk03} T. Blaszczyk, M. Ciesielski, M. Klimek and J. Leszczynski,
                    Numerical solution of fractional oscillator equation,
                    \emph{Appl. Math. Comput.} \textbf{218} (2011), 2480-2488.
\bibitem{Blaszczyk04} T. Blaszczyk and M. Ciesielski, Numerical solution of fractional
                    Euler-Lagrange equation with multipoint boundary conditions,
                     \emph{Scientific Research of the Institute of Mathematics and Computer 
                    Science} \textbf{2(10)} (2011), 43-48.
\bibitem{Blaszczyk05} T. Blaszczyk, J. Leszczynski and E. Szymanek, Numerical solution of 
                    composite left and right fractional Caputo derivative models 
                    for granular heat flow, \emph{Mech. Re. Commun.} 
                    \textbf{48} (2013), 42-45.
\bibitem{Diethelm01} K. Diethelm, {\it The Analysis of Fractional Differential Equations. 
                    An Application-Oriented Exposition Using Differential Operators of Caputo Type},
                    Springer-Verlag, Heidelberg etc  (2010).
\bibitem{Hilfer} R. Hilfer, {\it Applications of Fractional Calculus in
                    Physics}, World Scientific, Singapore (2000).
\bibitem{Kilbas}    A. A. Kilbas, H. M. Srivastava and J. J. Trujillo,
                    {\it Theory and Applications of Fractional Differential
                    Equations}, Elsevier, Amsterdam (2006).
\bibitem{Kiryakova} V. Kiryakova, {\it Generalized fractional calculus and applications},
                    Pitman Research Notes in Mathematics Series, 301,
                    Longman Sci. Tech., Harlow (1994).
\bibitem{Klimek01}  M. Klimek, Fractional sequential mechanics - models with
                    symmetric fractional derivative,
                    \emph{Czech. J. Phys.} \textbf{51} (2001), 1348-1354.
\bibitem{Klimek04}  M. Klimek, {\it On Solutions of Linear Fractional 
                    Differential Equations of a Variational Type},
                    The Publishing Office of the Czestochowa University of
                    Technology, Czestochowa (2009).
\bibitem{Klimek05}  M. Klimek, Existence - uniqueness result for a certain 
                    equation of motion in fractional mechanics,
                    \emph{Bull. Pol. Acad. Sci.: Technical Sciences} 
                    \textbf{58} (2010), 573-581.
\bibitem{Klimek07} M. Klimek and O.P. Agrawal, Fractional Sturm?Liouville problem,
                   \emph{Computers and Mathematics with Applications}
                   \textbf{66} (2013), 795-812.
\bibitem{Klimek06}  M. Klimek and M. Lupa, Reflection symmetric formulation of generalized
                    fractional variational calculus,
                    \emph{Fractional Calculus and Applied Analysis} \textbf{16} 
                    (2013), 243-261.
\bibitem{Lazo01}	 M. J. Lazo and D. F. M. Torres, 
                   The DuBois-Reymond fundamental lemma of the fractional
                   calculus of variations and an Euler-Lagrange equation 
                   involving only derivatives of Caputo,
                   \emph{J. Optim. Theory Appl.} \textbf{156}, 1 (2013), 56-67.
\bibitem{Leszczynski02} J. S. Leszczynski, {\it An Introduction to Fractional 
                    Mechanics}, The Publishing Office of the Czestochowa 
                    University of Technology, Czestochowa (2011).
\bibitem{Leszczynski03} J. S. Leszczynski and T. Blaszczyk, Modeling the transition 
                    between stable and unstable operation while emptying a silo, 
                    \emph{Granular Matter} \textbf{13} (2011), 429-438.
\bibitem{Lotfi}     A. Lotfi and S. A. Yousefi, A numerical technique for solving 
                    a class of fractional variational problems, \emph{Journal of 
                  Computational and Applied Mathematics} \textbf{237(1)} (2013), 633-643.
\bibitem{Magin}     R. L. Magin , \emph{Fractional Calculus in Bioengineering},
                    Begell House Inc., Redding (2006).
\bibitem{Malinowska01} A. B. Malinowska and D. F. M. Torres,
                     {\it Introduction to the Fractional Calculus of
                     Variations},
                     Imperial College Press, London (2012).
\bibitem{Malinowska02} A. B. Malinowska and D. F. M. Torres, 
                      Fractional calculus of variations for a combined Caputo
                      derivative, 
                      \emph{Fract. Calc. Appl. Anal.} \textbf{14}, 4 (2011), 
                      523--537.
\bibitem{Malinowska03} A. B. Malinowska and D. F. M. Torres,
                      Towards a combined fractional mechanics and quantization,
                      \emph{Fract. Calc. Appl. Anal.} \textbf{15}, 3 (2012), 
                      407--417.
\bibitem{Odzijewicz01} T. Odzijewicz, A. B. Malinowska and D. F. M. Torres,
                    Fractional variational calculus with classical and combined 
                    Caputo derivatives,
                    \emph{Nonlinear Anal.: TMA} \textbf{75} (2012), 1507-1515.
\bibitem{Odzijewicz02} T. Odzijewicz, A. B. Malinowska and D. F. M. Torres,
                    Green's theorem for generalized fractional derivatives, 
                    \emph{Fract. Calc. Appl. Anal.} \textbf{16}, 1 (2013), 
                    64--75.
\bibitem{Oldham}    K. B. Oldham and J. Spanier, {\it The Fractional Calculus:
                    Theory and Applications of Differentiation and Integration 
		to Arbitrary Order}, Academic Press, San Diego (1974). 
\bibitem{Podlubny}  I. Podlubny, {\it Fractional Differential Equations},
                    Academic Press, San Diego (1999).
\bibitem{Press} W. H. Press, S. A. Teukolsky, W. T. Vetterling and B. P. Flannery, 
					{\it Numerical Recipes: The Art of Scientific Computing (3rd ed.)}, 
					Cambridge University Press, New York (2007).
\bibitem{Pooseh01}  S. Pooseh, R. Almeida and D. F. M. Torres, 
                    Discrete direct methods in the fractional calculus of 
                    variations,
                    \emph{Comput. Math. Appl.} \textbf{66}, (2013), 668--676.
\bibitem{Pooseh02}  S. Pooseh, R. Almeida and D. F. M. Torres,
                    Numerical approximations of fractional derivatives with 
                    applications,
                    \emph{Asian Journal of Control} \textbf{15}, 3 (2013),
                    698--712.
\bibitem{Pooseh03}  S. Pooseh, R. Almeida and D. F. M. Torres,
                    A discrete time method to the first variation of fractional 
                    order variational functionals,
                    \emph{Cent. Eur. J. Phys.} (2013),  
                    DOI:10.2478/s11534-013-0250-0 (online first, to appear)
\bibitem{Riewe01}   F. Riewe, Nonconservative Lagrangian and Hamiltonian
                    mechanics, \emph{Phys. Rev. E} \textbf{53} (1996), 1890-1899.
\bibitem{Rivero}    M. Rivero, J. J. Trujillo and M. P. Velasco
                    A fractional approach to the Sturm-Liouville problem,
                    \emph{Cent. Eur. J. Phys.}
                    (2013) doi: 10.2478/s11534-013-0216-2.
\bibitem{Scalas}    E. Scalas, R. Gorenflo and F. Mainardi, Fractional calculus and
                    continuous time finance,
                    \emph{Physica A} \textbf{284} (2000), 376--384.
\bibitem{Wang}      D. Wang and A. Xiao, Fractional variational integrators for fractional
                    Euler-Lagrange equations with holonomic constraints,
                     \emph{Commun. Nonlinear Sci. Numer. Simulat.} \textbf{18} (2013), 905-914.


\end{thebibliography}
\end{document}